\documentclass{article}

\usepackage{amssymb}			
\usepackage{latexsym}			

\input{diagrams}
\diagramstyle[scriptlabels]

\newcommand{\hyph}{\mbox{-}}			
\newcommand{\diagspace}{\mbox{\hspace{2em}}}	
\newcommand{\cat}[1]{\mathcal{#1}}	
\newcommand{\scat}[1]{\mathbb{#1}}	
\newcommand{\fcat}[1]{\mathbf{#1}}	
\newcommand{\url}[1]{\mbox{\tt #1}}	

\newcommand{\blob}{\raisebox{0.1em}{\ensuremath{\scriptscriptstyle\bullet}}}
\newcommand{\of}{\raisebox{0.07em}{\ensuremath{\scriptstyle\circ}}}
\newcommand{\sof}{\raisebox{0.08ex}{\ensuremath{\scriptscriptstyle\circ}}}
\newcommand{\op}{\mathrm{op}}		
\newcommand{\id}{\mathit{id}}		
\newcommand{\ob}{\mathrm{ob}}		
\newcommand{\comp}{\mathit{comp}}	
\newcommand{\ids}{\mathit{ids}}		
\newcommand{\fc}{\fcat{fc}}		
\newcommand{\Bim}[1]{\fcat{Bim}(#1)}	
\newcommand{\Mon}[1]{\fcat{Mon}(#1)}	
\newcommand{\ehom}[3]{#1[#2,#3]}	
\newcommand{\eend}[2]{#1[#2]}		
\newcommand{\Set}{\fcat{Set}}		
\newcommand{\Ab}{\fcat{Ab}}		
\newcommand{\Cat}{\fcat{Cat}}		
\newcommand{\Multicat}{\fcat{Multicat}}	
\newcommand{\Graph}{\fcat{Graph}}	
\newcommand{\Action}{\fcat{Action}}	
\newcommand{\Span}{\fcat{Span}}		
\newcommand{\Bimod}{\fcat{Bimod}}	
\newcommand{\Prof}{\fcat{Prof}}		

\newtheorem{thm}{Theorem}
\newtheorem{defn}[thm]{Definition}
\newtheorem{prop}[thm]{Proposition}

\newarrow{Equals}=====			
\newcommand{\go}{\rTo}			
\newcommand{\goby}[1]{\rTo^{#1}}	
\newcommand{\goesto}{\,\longmapsto\,}	
\newcommand{\vslob}[3]			
	{\left.
	\begin{diagram}[height=1.5em]
	#1		\\
	\dTo>{\,#2}	\\
	#3		\\
	\end{diagram}
	\right.}

\newlength{\gwidth}	
\newlength{\gvert}	
\newlength{\gdrop}	
\newlength{\gbaredrop}	
\newlength{\goffset}	
\newlength{\gtemp}	
\newcommand{\present}[1]{%
\makebox[1\gwidth]{%
\rule[-1\gdrop]{0ex}{1\gvert}%
\raisebox{-1\gbaredrop}{#1}}}
\newcommand{\cinitdims}[2]{%
\setlength{\unitlength}{1em}%
\setlength{\goffset}{.35\unitlength}%
\setlength{\gwidth}{#1\unitlength}%
\setlength{\gvert}{#2\unitlength}%
\setlength{\gdrop}{.5\gvert}%
\addtolength{\gdrop}{-1\goffset}%
\setlength{\gbaredrop}{1\gdrop}%
\addtolength{\gvert}{.6\unitlength}%
\addtolength{\gdrop}{.3\unitlength}}	
\newcommand{\abovepic}[1]{%
\settoheight{\gtemp}{\ensuremath{#1}}%
\addtolength{\gvert}{1\gtemp}%
\settodepth{\gtemp}{\ensuremath{#1}}%
\addtolength{\gvert}{1\gtemp}}
\newcommand{\belowpic}[1]{%
\settoheight{\gtemp}{\ensuremath{#1}}%
\addtolength{\gvert}{1\gtemp}%
\addtolength{\gdrop}{1\gtemp}%
\settodepth{\gtemp}{\ensuremath{#1}}%
\addtolength{\gvert}{1\gtemp}%
\addtolength{\gdrop}{1\gtemp}}
\newcommand{\cell}[4]{\put(#1,#2){\makebox(0,0)[#3]{\ensuremath{#4}}}}
\newcommand{\prectwo}[3]%
{\begin{picture}(4.2,3.4)(-0.1,-0.2)%
\cell{2}{3.2}{b}{#1}%
\cell{2}{-0.2}{t}{#2}%
\cell{2.2}{1.5}{l}{#3}%
\qbezier(0,2)(2,4)(4,2)%
\qbezier(0,1)(2,-1)(4,1)%
\put(4,2){\vector(1,-1){0}}%
\put(4,1){\vector(1,1){0}}%
\cell{2}{1.5}{c}{\Downarrow}%
\end{picture}}
\newcommand{\ctwo}[3]{%
\cinitdims{4.2}{3.4}%
\abovepic{#1}%
\belowpic{#2}%
\present{\prectwo{#1}{#2}{#3}}}

\begin{document}

\title{Generalized Enrichment of Categories}
\author{Tom Leinster\thanks{Supported by the EPSRC and St John's College,
Cambridge}\\ \\
	\normalsize{Department of Pure Mathematics and Mathematical Statistics,}\\
	\normalsize{Centre for Mathematical Sciences, Wilberforce Road,}\\ 
	\normalsize{Cambridge CB3 0WB, UK}\\
	\normalsize{\url{leinster@dpmms.cam.ac.uk}}}
\date{}

\maketitle

\vspace*{-1em}
\begin{center}
\textbf{Abstract}
\end{center}
\begin{quote}
We define the phrase `category enriched in an $\fc$-multicategory' and explore
some examples. An $\fc$-multicategory is a very general kind of 2-dimensional
structure, special cases of which are double categories, bicategories,
monoidal categories and ordinary multicategories. Enrichment in an
$\fc$-multicategory extends the (more or less well-known) theories of
enrichment in a monoidal category, in a bicategory, and in a
multicategory. Moreover, $\fc$-multicategories provide a natural setting for
the bimodules construction, traditionally performed on suitably cocomplete
bicategories. Although this paper is elementary and self-contained, we also
explain why, from one point of view, $\fc$-multicategories are
the \emph{natural} structures in which to enrich categories.

\small
Published as \emph{Journal of Pure and Applied Algebra} \textbf{168}
(2002), 391--406. 

\end{quote}

A general question in category theory is: given some kind of categorical
structure, what might it be enriched in? For instance, suppose we take
braided monoidal categories. Then the question asks: what kind of thing must
$\cat{V}$ be if we are to speak sensibly of $\cat{V}$-enriched braided monoidal
categories? (The usual answer is that $\cat{V}$ must be a symmetric monoidal
category.) 

In another paper, \cite{GECM}, I have given an answer to the general question
for a certain family of categorical structures (generalized
multicategories). In particular, this theory gives an answer to the question
`what kind of structure $\cat{V}$ can a category be enriched in'? The answer
is: an `$\fc$-multicategory'.

Of course, the traditional answer to this question is that $\cat{V}$ is a
monoidal category. But there is also a notion of a category enriched in a
bicategory (see Walters \cite{Wal}). And generalizing in a different
direction, it is easy to see how one might speak of a category enriched in an
ordinary multicategory (`change tensors to commas'). An $\fc$-multicategory is,
in fact, a very general kind of 2-dimensional categorical structure,
encompassing monoidal categories, bicategories, multicategories and double
categories. The theory of categories enriched in an $\fc$-multicategory extends
all of the aforementioned theories of enrichment.

So from the point of view of \cite{GECM}, $\fc$-multicategories are the natural
structures in which to enrich a category. In this work, however, we do not
assume any knowledge of \cite{GECM} or of generalized
multicategories. Instead, we define $\fc$-multicategory in an elementary
fashion (Section~\ref{sec:fcms}) and then define what a category enriched in
an $\fc$-multicategory is (Section~\ref{sec:enr}). Along the way we see how
enrichment in an $\fc$-multicategory extends the previously-mentioned theories
of enrichment, and look at various examples.

$\fc$-multicategories also provide a natural setting for the bimodules
construction (Section~\ref{sec:bim}), traditionally carried out on
bicategories satisfying certain cocompleteness conditions. At the level of
$\fc$-multicategories, the construction is both more general and free of
technical restrictions. We show, in particular, that a category enriched in
an $\fc$-multicategory $\cat{V}$ naturally gives rise to a category enriched
in the $\fc$-multicategory $\Bim{\cat{V}}$ of bimodules in $\cat{V}$. This
result is functorial (that is, a $\cat{V}$-enriched functor gives rise to a
$\Bim{\cat{V}}$-enriched functor), a statement which only holds if we work
with $\fc$-multicategories rather than bicategories.

\section{$\fc$-multicategories}	\label{sec:fcms}

In a moment, an explicit and elementary definition of $\fc$-multicategory will
be given. But first it might be helpful to look briefly at the wider
context in which this definition sits: the theory of `generalized
multicategories'. The reader is reassured that no knowledge of this wider
context is required in order to understand the rest of the paper.

Given a monad $T$ on a category $\cat{E}$, both having certain convenient
properties, there is a category of \emph{$T$-multicategories}. A
$T$-multicategory $C$ consists of a diagram
\begin{diagram}[size=2em]
	&		&C_1	&		&	\\
	&\ldTo<{\mathrm{dom}}&	&\rdTo>{\mathrm{cod}}&	\\
T(C_0)	&		&	&		&C_0	\\
\end{diagram}
in $\cat{E}$ (a \emph{$T$-graph}) together with functions defining
`composition' and `identity'; the full details can be found in
Burroni~\cite{Bur} or Leinster (\cite{GOM} or \cite{OHDCT}). Thus when $T$ is
the identity monad on $\cat{E}=\Set$, a $T$-multicategory is simply a
category. When $T$ is the free-monoid monad on $\cat{E}=\Set$, a
$T$-multicategory is a multicategory in the original sense of
Lambek~\cite{Lam}. When $T$ is the free (strict) $\infty$-category monad on
the category $\cat{E}$ of globular sets (`$\infty$-graphs'), a
$T$-multicategory $C$ with $C_0=1$ is a higher operad in the sense of
Batanin~\cite{Bat}. The example which concerns us here is when $T$ is the
free category monad $\fc$ on the category $\cat{E}$ of directed graphs. A
$T$-multicategory is then an $\fc$-multicategory in the sense of the
following explicit definition.

\begin{defn}	\label{defn:fcm}
An \emph{$\fc$-multicategory} consists of
\begin{itemize}
\item A class of \emph{objects} $x$, $x'$, \ldots
\item For each pair $(x,x')$ of objects, a class of \emph{vertical
1-cells} $\vslob{x}{}{x'}$, denoted $f$, $f'$, \ldots
\item For each pair $(x,x')$ of objects, a class of \emph{horizontal
1-cells} $x\go x'$, denoted $m$, $m'$, \ldots
\item For each $n\geq 0$, objects $x_0,\ldots,x_n,x,x'$, vertical 1-cells
$f,f'$, and horizontal 1-cells $m_1,\ldots,m_n,m$, a class of
\emph{2-cells}
\begin{equation}	\label{eq:two-cell}
\begin{diagram}[height=2em]
x_0	&\rTo^{m_1}	&x_1	&\rTo^{m_2}	&\ 	&\cdots	
&\ 	&\rTo^{m_n}	&x_n	\\
\dTo<{f}&		&	&		&\Downarrow	&
&	&		&\dTo>{f'}\\
x	&		&	&		&\rTo_{m}	&	
&	&		&x',	\\
\end{diagram}
\end{equation}
denoted $\theta$, $\theta'$, \ldots
\item \emph{Composition} and \emph{identity} functions making the objects
and vertical 1-cells into a category
\item A \emph{composition} function for 2-cells, as in the picture
\begin{diagram}[width=.5em,height=1em]
\blob&\rTo^{m_1^1}&\cdots&\rTo^{m_1^{r_1}}&
\blob&\rTo^{m_2^1}&\cdots&\rTo^{m_2^{r_2}}&\blob&
\ &\cdots&\ &
\blob&\rTo^{m_n^1}&\cdots&\rTo^{m_n^{r_n}}&\blob\\
\dTo<{f_0}&&\Downarrow\,\scriptstyle{\theta_1}&&
\dTo&&\Downarrow\,\scriptstyle{\theta_2}&&\dTo&
\ &\cdots&\ &
\dTo&&\Downarrow\,\scriptstyle{\theta_n}&&\dTo>{f_n}\\
\blob&&\rTo_{m_1}&&
\blob&&\rTo_{m_2}&&\blob&
\ &\cdots&\ &
\blob&&\rTo_{m_n}&&\blob\\
\dTo<{f}&&&&&&&&\Downarrow\,\scriptstyle{\theta} &&&&&&&&\dTo>{f'}\\
\blob&&&&&&&&\rTo_{m}&&&&&&&&\blob\\
\end{diagram}
$\goesto$
\begin{diagram}[width=.5em,height=1em]
\blob&\rTo^{m_1^1}&\ &&
&&&&\cdots&
&&&
&&\ &\rTo^{m_n^{r_n}}&\blob\\
&&&&&&&&&&&&&&&&\\
\dTo<{f\of f_0}&&&&&&&&\Downarrow\,\scriptstyle{\theta\of(\theta_1,\theta_2,\ldots,\theta_n)}
&&&&&&&&\dTo>{f'\of f_n}\\
&&&&&&&&&&&&&&&&\\
\blob&&&&&&&&\rTo_{m}&&&&&&&&\blob\\
\end{diagram}
($n\geq 0, r_i\geq 0$), where the $\blob$'s represent objects 
\item An \emph{identity} function
\[
\begin{diagram}[width=1em,height=2em]
x&\rTo^{m}&x'\\
\end{diagram}
\diagspace\goesto\diagspace
\begin{diagram}[width=1em,height=2em]
x&\rTo^{m}&x'\\
\dTo<{1_x}&\Downarrow\,\scriptstyle{1_{m}}&\dTo>{1_{x'}}\\
x&\rTo_{m}&x'.\\
\end{diagram}
\]
\end{itemize}
The 2-cell composition and identities are required to obey associativity and
identity laws.
\end{defn}
The associativity and identity laws ensure that any diagram of
pasted-together 2-cells with a rectangular boundary has a well-defined
composite.

\subsection*{Examples}

\begin{enumerate}

\item		\label{eg:double}
Any double category gives an $\fc$-multicategory, in which a 2-cell as in
diagram~(\ref{eq:two-cell}) is a 2-cell
\begin{diagram}[height=2em]
x_0	&\rTo^{m_n \sof \cdots \sof m_1}&x_n\\
\dTo<{f}&\Downarrow			&\dTo>{f'}\\
x	&\rTo_{m}			&x'\\
\end{diagram}
in the double category. If the double category is called $\cat{D}$ then we also
call the resulting $\fc$-multicategory $\cat{D}$, and we use the same
convention for bicategories (next example).

\item	\label{eg:bicat}
Any bicategory gives an $\fc$-multicategory in which the only vertical 1-cells
are identity maps, and a 2-cell as in diagram~(\ref{eq:two-cell}) is a 2-cell
\[
x_0
\ctwo{\scriptstyle{m_n \sof\cdots\sof m_1}}{\scriptstyle{m}}{}
x_n
\]
in the bicategory (with $x=x_0$ and $x'=x_n$). 

Here $m_n \of\cdots\of m_1$ denotes some $n$-fold composite of the 1-cells
$m_n, \ldots, m_1$ in the bicategory. For the sake of argument let us decide
to associate to the left, so that $m_4\of m_3\of m_2\of m_1$ means $((m_4\of
m_3)\of m_2)\of m_1$. A different choice of bracketing would only affect the
resulting $\fc$-multicategory up to isomorphism (in the obvious sense).

\item	\label{eg:mon-cat}
Any monoidal category $\cat{M}$ gives rise to an $\fc$-multicategory
$\Sigma\cat{M}$ (the \emph{suspension} of $\cat{M}$) in which there is one
object and one vertical 1-cell, and a 2-cell
\begin{equation}	\label{eq:vt-two-cell}
\begin{diagram}[size=2em,abut]
\bullet	&\rLine^{M_1}	&\bullet	&\rLine^{M_2}	&\bullet	&\cdots 
&\bullet	&\rLine^{M_n}	&\bullet	\\
\dLine<1&		&	&		&\Downarrow&
&	&		&\dLine>1\\
\bullet	&		&	&		&\rLine_{M}&
&	&		&\bullet\\
\end{diagram}
\end{equation}
is a morphism $M_n\otimes\cdots\otimes M_1 \go M$ in $\cat{M}$. This is a
special case of Example~(\ref{eg:bicat}).

\item	\label{eg:multicat}
Similarly, any ordinary multicategory $\cat{M}$ gives an $\fc$-multicategory
$\Sigma\cat{M}$: there is one object, one vertical 1-cell, and a 2-cell as in
diagram~(\ref{eq:vt-two-cell}) is a map $M_1,\ldots,M_n \go M$ in $\cat{M}$.

\item \label{eg:Span}
We define an $\fc$-multicategory $\Span$. Objects are sets, vertical 1-cells
are functions, a horizontal 1-cell $X\go Y$ is a diagram
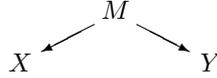
\[
\begin{diagram}[width=1.8em,height=1em]
	&	&M	&	&	\\
	&\ldTo	&	&\rdTo	&	\\
X	&	&	&	&Y	\\
\end{diagram}
\]
of sets and functions, and a 2-cell inside
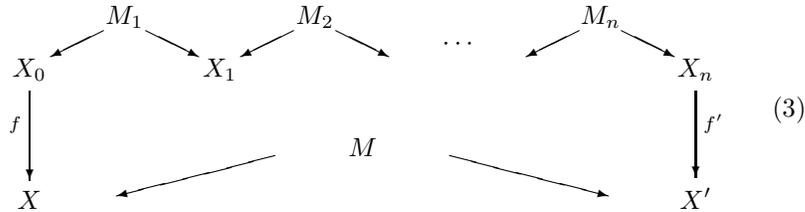
\begin{equation}	\label{eq:span-two-cell}
\begin{diagram}[width=1.8em,height=1em,tight]
	&	&M_1	&	&	&	&M_2	&	&	
&	&	&	&M_n	&	&	\\
	&\ldTo	&	&\rdTo	&	&\ldTo	&	&\rdTo	&
&\cdots	&	&\ldTo	&	&\rdTo	&	\\ 
X_0	&	&	&	&X_1	&	&	&	&\ 
&	&\ 	&	&	&	&X_n	\\
	&	&	&	&	&	&	&	&
&	&	&	&	&	&	\\
\dTo<f	&	&	&	&	&	&	&	&
&	&	&	&	&	&\dTo>{f'}\\
	&	&	&	&	&	&	&M	&
&	&	&	&	&	&	\\
	&	&	&	&	&	&\ldTo(7,2)&	&\rdTo(7,2)
&	&	&	&	&	&	\\
X	&	&	&	&	&	&	&	&
&	&	&	&	&	&X'	\\
\end{diagram}
\end{equation}
is a function $\theta$ making
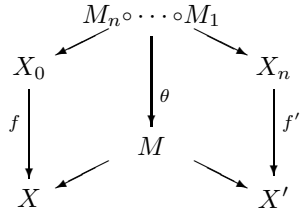
\begin{diagram}[width=2em,height=1em]
	&	&M_n\of\cdots\of M_1	&	&	\\
	&\ldTo	&			&\rdTo	&	\\
X_0	&	&			&	&X_n	\\
	&	&\dTo>{\theta}		&	&	\\
\dTo<f	&	&			&	&\dTo>{f'}\\
	&	&M			&	&	\\
	&\ldTo	&			&\rdTo	&	\\
X	&	&			&	&X'	\\
\end{diagram}
commute. Here $M_n\of\cdots\of M_1$ is the limit of the top row of
diagram~(\ref{eq:span-two-cell}), an iterated pullback. Composition is
defined in the obvious way.

$\Span$ is an example of a `weak double category', which is just like a
double category except that horizontal 1-cell composition only obeys
associativity and identity axioms up to coherent isomorphism.

It is rather idiosyncratic to name this $\fc$-multicategory after its
horizontal 1-cells: usually one names a categorical structure after its
objects (e.g.\ $\fcat{Group}$, $\Set$). However, we do not want to confuse
the $\fc$-multicategory $\Span$ of sets with the mere category $\Set$ of
sets, so we will stick to this convention.

Notice, incidentally, that $\Set$ is the category formed by the objects and
vertical 1-cells of $\Span$, and that the $\fc$-multicategory $\Sigma\Set$
arising from the monoidal category $(\Set,\times,1)$ is the `full'
sub-$\fc$-multicategory of $\Span$ whose only object is $1$.

\item	\label{eg:Prof}
There is an $\fc$-multicategory $\Prof$, in which the category formed by the
objects and vertical 1-cells is the usual category of (small) categories and
functors. Horizontal 1-cells are profunctors (bimodules): that is, a
horizontal 1-cell $\scat{X} \go \scat{X'}$ is a functor $\scat{X}^\op \times
\scat{X'} \go \Set$. A 2-cell
\begin{diagram}[height=2em]
\scat{X}_0	&\rTo^{M_1}	&\scat{X}_1	&\rTo^{M_2}	&\ 	&\cdots	
&\ 	&\rTo^{M_n}	&\scat{X}_n	\\
\dTo<{F}&		&	&		&\Downarrow	&
&	&		&\dTo>{F'}\\
\scat{X}	&		&	&		&\rTo_{M}	&	
&	&		&\scat{X}'	\\
\end{diagram}
consists of a function
\[
M_n(x_{n-1}, x_n) \times\cdots\times M_1(x_0, x_1)
\go
M(F(x_0), F'(x_n))
\]
for each $x_0\in\scat{X}_0, \ldots, x_n\in\scat{X}_n$, such that this family
of functions is natural in the $x_i$'s. So if the functors $F$ and $F'$ are
identities then this is a morphism of profunctors $M_n \otimes \cdots \otimes
M_1 \go M$.

\item	\label{eg:Bimod}
In a similar spirit, $\Bimod$ is the following $\fc$-multicategory:
\begin{itemize}
\item
objects are rings (with identity, not necessarily commutative)
\item
vertical 1-cells are ring homomorphisms
\item
a horizontal 1-cell $R\go S$ is an $(S,R)$-bimodule
\item
a 2-cell
\begin{diagram}[height=2em]
R_0	&\rTo^{M_1}	&R_1	&\rTo^{M_2}	&\ 	&\cdots	
&\ 	&\rTo^{M_n}	&R_n	\\
\dTo<{f}&		&	&		&\Downarrow\,\scriptstyle{\theta}&
&	&		&\dTo>{f'}\\
R	&		&	&		&\rTo_{M}	&	
&	&		&R'	\\
\end{diagram}
is a function $\theta:M_n\times\cdots\times M_1 \go M$ which preserves
addition in each component separately (is `multi-additive') and satisfies the
equations 
\begin{eqnarray*}
\theta(r_n\cdot m_n, m_{n-1}, \ldots)		&=&
f'(r_n)\cdot \theta(m_n, m_{n-1}, \ldots)	\\
\theta(m_n\cdot r_{n-1}, m_{n-1}, \ldots)	&=&
\theta(m_n, r_{n-1}\cdot m_{n-1}, \ldots)	\\
	&\vdots		\\
\theta(\ldots, m_2\cdot r_1, m_1)		&=&
\theta(\ldots, m_2, r_1\cdot m_1)		\\
\theta(\ldots, m_2, m_1\cdot r_0)		&=&
\theta(\ldots, m_2, m_1)\cdot f(r_0)
\end{eqnarray*} 
\item
composition and identities are defined in the evident way.
\end{itemize}

\item	\label{eg:Action}
If we remove all the additive structure involved in $\Bimod$ then we
obtain an $\fc$-multicategory $\Action$; alternatively, $\Action$ is the
`full' sub-$\fc$-multicategory of $\Prof$ in which the only objects allowed
are 1-object categories. Thus the objects of $\Action$ are monoids, the
vertical 1-cells are monoid homomorphisms, a horizontal 1-cell $R\go S$ is a
set with commuting left $S$-action and right $R$-action, and 2-cells are
defined as in Example~(\ref{eg:Bimod}).

\end{enumerate}

\section{Enrichment}	\label{sec:enr}

The purpose of this paper is to explore in an elementary way the concept of a
category enriched in an $\fc$-multicategory. But just as the elementary
definition of $\fc$-multicategory (Definition~\ref{defn:fcm}) is plucked out of
a much larger theory (as explained in the introduction to
Section~\ref{sec:fcms}), so too is the definition of category enriched in an
$\fc$-multicategory. There is a whole theory \cite{GECM} of enrichment for
generalized multicategories, of which the present work is just the most
simple case. This wider theory runs as follows.

Any $T$-multicategory has an underlying $T$-graph, as explained above, and so
there is a forgetful functor
\[
T\hyph\Multicat \go T\hyph\Graph.
\]
Under mild conditions on $\cat{E}$ and $T$, this functor has a left
adjoint. We thus obtain a monad $T'$ on the category
$\cat{E'}=T\hyph\Graph$. We can then speak of $T'$-multicategories, and if
$\cat{V}$ is a $T'$-multicategory one can make a definition of
\emph{$\cat{V}$-enriched $T$-multicategory}. So: we can speak of a
$T$-multicategory enriched in a $T'$-multicategory.

The most simple case is the identity monad $T$ on $\cat{E}=\Set$. Then
$T$-multicategories are categories, $T'$ is the free category monad $\fc$ on
$\cat{E'}=\Graph$, and $T'$-multicategories are $\fc$-multicategories. So the
general theory gives a concept of category enriched in an $\fc$-multicategory.
The main part of this section is a direct description of this concept.

The next most simple case is the free monoid monad $T$ on $\cat{E}=\Set$, and
here there are two especially interesting examples of enriched
$T$-multicategories. Firstly, it turns out that any symmetric monoidal
category $\cat{S}$ gives rise to a $T'$-multicategory $\cat{V}$, and a
one-object $\cat{V}$-enriched $T$-multicategory is then exactly what
topologists call a (non-symmetric) operad in $\cat{S}$ (see e.g.\
\cite{MayDOA}). Secondly, there is a certain naturally-arising
$T'$-multicategory $\cat{V}$ such that $\cat{V}$-enriched $T$-multicategories
are the structures called `relaxed multicategories' by Borcherds in his
definition of vertex algebras over a vertex group (\cite{Bor}, \cite{SnyEBG},
\cite{SnyRMS}), and called `pseudo-monoidal categories' by Soibelman in his
work on quantum affine algebras (\cite{SoiMTC}, \cite{SoiMBC}).

The general definition of enriched $T$-multicategory is very simple. Take a
monad $T$ on a category $\cat{E}$, and let $T'$ be the free $T$-multicategory
monad, as above. Given an object $A$ of $\cat{E}$, we can form $I(A)$ (with $I$
for `indiscrete'), the unique $T$-multicategory with graph
\[
T(A) \lTo^{\mathrm{pr}_1} T(A) \times A \rTo^{\mathrm{pr}_2} A.
\]
Then $I(A)$ is a $T'$-algebra, say $h: T'(I(A)) \go I(A)$. Arising from this
is a $T'$-multicategory $M(I(A))$, the unique such with graph
\[
T'(I(A)) \lTo^{1} T'(I(A)) \rTo^{h} I(A).
\]
For a fixed $T'$-multicategory $\cat{V}$, a \emph{$\cat{V}$-enriched
$T$-multicategory} is defined as an object $C_0$ of $\cat{E}$ together with a
map $T'(I(C_0)) \go \cat{V}$ of $T'$-multicategories. Maps between
$\cat{V}$-enriched $T$-multicategories are also defined in a simple way, thus
giving a category.

In the case concerning us, $\cat{E}=\Set$ and $T=\id$, the definition of
enriched ($T$-multi)category is therefore as follows. Given a set $A$, we
obtain the indiscrete category $I(A)$ on $A$. In the $\fc$-multicategory
$M(I(A))$, an object is an element of $A$, the only vertical 1-cells are
identities, there is one horizontal 1-cell $a\go b$ for each $a,b\in A$, and
for each $a_0, \ldots, a_n \in A$ there is precisely one 2-cell of the form 
\[
\begin{diagram}[height=2em]
a_0	&\rTo		&a_1	&\rTo		&\ 	&\cdots	
&\ 	&\rTo		&a_n	\\
\dTo<{1}&		&	&		&\Downarrow	&
&	&		&\dTo>{1}\\
a_0	&		&	&		&\rTo		&	
&	&		&a_n.	\\
\end{diagram}
\]
Composition and identities are uniquely determined. A category enriched in an
$\fc$-multicategory $\cat{V}$ consists of a set $C_0$ together with a map
from the $\fc$-multicategory $M(I(C_0))$ to $\cat{V}$. This definition is
plainly equivalent to Definition~\ref{defn:enr-cat} below.

That concludes the sketch of the theory of enriched generalized
multicategories, and we now return to the elementary account.

Fix an $\fc$-multicategory $\cat{V}$.

\begin{defn}	\label{defn:enr-cat}
A \emph{category enriched in $\cat{V}$}, or \emph{$\cat{V}$-enriched
category}, $C$, consists of
\begin{itemize}
\item 
a class $C_0$ (of `objects')
\item 
for each $a\in C_0$, an object $\eend{C}{a}$ of $\cat{V}$
\item 
for each $a,b\in C_0$, a horizontal 1-cell $\eend{C}{a}
\rTo^{\ehom{C}{a}{b}} \eend{C}{b}$ in $\cat{V}$
\item 
for each $a,b,c\in C_0$, a `composition' 2-cell
\begin{diagram}[height=2em]
\eend{C}{a}&\rTo^{\ehom{C}{a}{b}}&\eend{C}{b}	&\rTo^{\ehom{C}{b}{c}}
&\eend{C}{c}	\\
\dTo<{1}&		&\Downarrow\,\scriptstyle{\comp_{a,b,c}}	&	&\dTo>{1}\\
\eend{C}{a}&		&\rTo_{\ehom{C}{a}{c}}		&	&\eend{C}{c}\\
\end{diagram}
\item
for each $a\in C_0$, an `identity' 2-cell
\begin{diagram}[height=2em]
\eend{C}{a}	&\rEquals		&\eend{C}{a}	\\
\dTo<{1}	&\Downarrow\,\scriptstyle{\ids_{a}}	&\dTo>{1}	\\
\eend{C}{a}	&\rTo_{\ehom{C}{a}{a}}	&\eend{C}{a}	\\
\end{diagram}
(where the equality sign along the top denotes a string of 0 horizontal
1-cells)
\end{itemize}
such that $\comp$ and $\ids$ satisfy associativity and identity axioms.
\end{defn}

To the reader used to enrichment in a monoidal category, the only unfamiliar
piece of data in this definition is the family of objects $\eend{C}{a}$. To
the reader used to enrichment in bicategories even this will be familiar;
indeed, since the vertical 1-cells are not used in any significant way, our
definition looks very much like the definition of category enriched in a
bicategory (see \cite{Wal}). This lack of use of the vertical 1-cells might
seem to weigh against the claim that $\fc$-multicategories are, in some
sense, the natural structures in which to enrich categories. However, the
vertical 1-cells \emph{are} used in the definition of $\cat{V}$-enriched
functor, which is given next. This makes the theory of enrichment in an
$\fc$-multicategory run more smoothly (sometimes, at least) than that of
enrichment in a bicategory, as we shall see towards the end of
Section~\ref{sec:bim}.

\begin{defn}	\label{defn:enr-ftr}
Let $C$ and $D$ be $\cat{V}$-enriched categories. A \emph{$\cat{V}$-enriched
functor} $F:C\go D$ consists of 
\begin{itemize}
\item
a function $F:C_0 \go D_0$
\item
for each $a\in C_0$, a vertical 1-cell
$\vslob{\eend{C}{a}}{F_a}{\eend{D}{F(a)}}$ 
\item
for each $a,b\in C_0$, a 2-cell
\begin{diagram}[height=2em]
\eend{C}{a}	&\rTo^{\ehom{C}{a}{b}}		&\eend{C}{b}	\\
\dTo<{F_a}	&\Downarrow\,\scriptstyle{F_{ab}}&\dTo>{F_b}	\\
\eend{D}{F(a)}	&\rTo_{\ehom{D}{F(a)}{F(b)}}	&\eend{D}{F(b)}
\end{diagram}
\end{itemize}
such that the $F_{ab}$'s commute with the composition and identity 2-cells in
$C$ and $D$, in an evident sense. 
\end{defn}

With the obvious notion of composition of $\cat{V}$-enriched functors, we
obtain a category $\cat{V}\hyph\Cat$ of $\cat{V}$-enriched categories and
functors.

\subsection*{Examples}

\begin{enumerate}

\item
Let $\cat{M}$ be a monoidal category and consider a category $C$ enriched in
the $\fc$-multicategory $\Sigma\cat{M}$ (defined in
Example~\ref{sec:fcms}(\ref{eg:mon-cat})). There is only one possible choice
for the $\eend{C}{a}$'s, so the data for $C$ consists of the set $C_0$, the
objects $\ehom{C}{a}{b}$ of $\cat{M}$, and the maps
\[
\ehom{C}{b}{c} \otimes \ehom{C}{a}{b} \go \ehom{C}{a}{c},
\diagspace
I\go\ehom{C}{a}{a}.
\]
Thus it turns out that a category enriched in $\Sigma\cat{M}$ is just a
category enriched (in the usual sense) in $\cat{M}$. The same goes for
enriched functors, so $(\Sigma\cat{M})\hyph\Cat$ is isomorphic to the usual
category of $\cat{M}$-enriched categories and functors.

\item
Let $\cat{M}$ be an (ordinary) multicategory. There is an obvious notion of
category enriched in $\cat{M}$: that is, a set $C_0$ together with an object
$\ehom{C}{a}{b}$ of $\cat{M}$ for each $a,b\in C_0$ and arrows
\[
\ehom{C}{a}{b}, \ehom{C}{b}{c} \go \ehom{C}{a}{c},
\diagspace
\cdot \go \ehom{C}{a}{a}
\]
(where $\cdot$ is the empty sequence), obeying suitable axioms. This is
precisely the same thing as a category enriched in $\Sigma\cat{M}$.

\item
If $\cat{B}$ is a bicategory then our $\cat{B}\hyph\Cat$ is isomorphic to the
category of $\cat{B}$-enriched categories defined in Walters \cite{Wal}.

\item
Fix a topological space $A$. Then there is a bicategory $\Pi_2 A$, the
\emph{homotopy bicategory} of $A$, in which an object is a point of $A$, a
1-cell is a path in $A$, and a 2-cell is a homotopy class of path homotopies
in $A$. For any 1-cell $\gamma:a \go b$ there is an associated 1-cell
$\gamma^*: b\go a$ (that is, $\gamma$ run backwards), and there are canonical
2-cells $1_b \go \gamma\of\gamma^*$ and $\gamma^*\of \gamma \go 1_a$. 

Now suppose that $A$ is nonempty and path-connected, and make a choice of a
basepoint $a_0$ and for each $a\in A$ a path $\gamma_a: a_0 \go a$. Then we
obtain a category $C$ enriched in $\Pi_2 A$, as follows:
\begin{itemize}
\item
$C_0$ is the underlying set of $A$
\item
$\eend{C}{a}=a$
\item
$\ehom{C}{a}{b} = \gamma_b \of \gamma_a^*$ (a path from $a$ to $b$)
\item
composition $\ehom{C}{b}{c}\of\ehom{C}{a}{b} \go \ehom{C}{a}{c}$ is the
2-cell 
\[
(\gamma_c \of \gamma_b^*) \of (\gamma_b \of \gamma_a^*) 
\go \gamma_c \of \gamma_a^*
\]
 coming from the canonical 2-cell $\gamma_b^* \of \gamma_b \go
1_{a_0}$ 
\item
the identity 2-cell $1_a\go \ehom{C}{a}{a}$ is the canonical 2-cell $1_a \go
\gamma_a \of \gamma_a^*$.
\end{itemize}

\item
In the previous example, the bicategory $\Pi_2 A$ can be replaced by any
bicategory $\cat{B}$ in which the underlying directed graph of objects and
1-cells is (nonempty and) connected, and every 1-cell has a left
adjoint. (I thank the referee for alerting me to this.) 

\item	\label{eg:comma}
$\Span\hyph\Cat$ is equivalent to the comma category $(\ob\downarrow\Set)$,
where $\ob:\Cat\go\Set$ is the objects functor. This means that a category
enriched in $\Span$ consists of a category $D$, a set $I$, and a function
$\ob(D) \go I$. To see why this is true, recall that a category $C$ enriched
in $\Span$ consists of
\begin{itemize}
\item
a set $C_0$
\item
for each $i\in C_0$, a set $\eend{C}{i}$
\item
for each $i,j\in C_0$, a span 
\begin{diagram}
\eend{C}{i} &\lTo^{s_{ij}} &\ehom{C}{i}{j} &\rTo^{t_{ij}} &\eend{C}{j}
\end{diagram}
\item
composition functions $\ehom{C}{j}{k} \times_{\eend{C}{j}} \ehom{C}{i}{j} \go
\ehom{C}{i}{k}$
\item
identity functions $\eend{C}{i} \go \ehom{C}{i}{i}$,
\end{itemize}
all satisfying axioms. We can construct from $C$ a category $D$ with
object-set $\coprod_{i\in C_0}\eend{C}{i}$, arrow-set $\coprod_{i,j\in
C_0}\ehom{C}{i}{j}$, source and target maps given by the $s_{ij}$'s and
$t_{ij}$'s, and composition and identity operations coming from those in
$C$. By taking $I=C_0$ and the projection function $\ob(D) \go I$, we now
have an object of $(\ob\downarrow\Set)$. A similar analysis of
$\Span$-enriched functors can be carried out, and we end up with a functor
\[
\Span\hyph\Cat \go (\ob \downarrow \Set).
\]
It is easy to see that this functor is an equivalence.

\end{enumerate}

Let us briefly consider enriched categories with only one object. In the
classical case of enrichment in a monoidal category $\cat{M}$, the category
of one-object $\cat{M}$-enriched categories is the category $\Mon{\cat{M}}$
of monoids in $\cat{M}$. For an arbitrary $\fc$-multicategory $\cat{V}$, we
therefore define $\Mon{\cat{V}}$ to be the full subcategory of
$\cat{V}\hyph\Cat$ whose objects are $\cat{V}$-enriched categories $C$ with
$|C_0| = 1$. Definitions~\ref{defn:enr-cat} and~\ref{defn:enr-ftr} yield an
explicit description of $\Mon{\cat{V}}$.

\subsection*{Examples}

\begin{enumerate}

\item
If $\cat{M}$ is a monoidal category then $\Mon{\Sigma\cat{M}}$ is the
category of monoids in $\cat{M}$.

\item
If $\cat{M}$ is a multicategory then an object of $\Mon{\Sigma\cat{M}}$
consists of an object $M$ of $\cat{M}$ together with maps
\[
M,M\go M, \diagspace \cdot\go M
\]
satisfying associativity and identity laws---in other words, a `monoid in
$\cat{M}$'. A monoid in $\cat{M}$ is also the same thing as a multicategory map
$\fcat{1}\go\cat{M}$, where $\fcat{1}$ is the terminal multicategory. 

\item 	\label{eg:strict-monads} 
If $\cat{B}$ is a bicategory then an object of $\Mon{\cat{B}}$ is a monad in
$\cat{B}$ in the sense of Street \cite{St}: that is, it's an object $X$ of
$\cat{B}$ together with a 1-cell $t:X\go X$ and 2-cells $\mu: t\of t \go t$,
$\eta: 1 \go t$ satisfying the usual monad axioms. There are no maps
$(X,t,\mu,\eta) \go (X',t',\mu',\eta')$ in $\Mon{\cat{B}}$ unless $X=X'$, and
in this case such a map consists of a 2-cell $t\go t'$ commuting with the
$\mu$'s and $\eta$'s. So $\Mon{\cat{B}}$ is the category of monads and
`strict monad maps' in $\cat{B}$.

\item
Let $\cat{B}$ be a 2-category. Associated to $\cat{B}$ is not only the
$\fc$-multicategory $\cat{B}$ of the previous example---which we now call
$\cat{V}$---but also two more $\fc$-multicategories, $\cat{W}$ and
$\cat{W'}$. Both $\cat{W}$ and $\cat{W'}$ are defined from double categories
(see Example~\ref{sec:fcms}(\ref{eg:double})), and in both cases an object is
an object of $\cat{B}$, a vertical 1-cell is a 1-cell of $\cat{B}$, and a
horizontal 1-cell is also a 1-cell of $\cat{B}$. In the case of $\cat{W}$, a
2-cell inside
\begin{diagram}[size=2em]
X	&\rTo^{t}	&Y		\\
\dTo<f	&		&\dTo>{g}	\\
X'	&\rTo_{t'}	&Y'		
\end{diagram}
is a 2-cell $t'\of f \go g\of t$ in $\cat{B}$; in the case of $\cat{W'}$, it
is a 2-cell $g\of t \go t'\of f$ in $\cat{B}$. Composition and identities are
defined in the obvious way.

Since $\cat{V}$, $\cat{W}$ and $\cat{W'}$ are identical when we ignore the
vertical 1-cells, the objects of $\Mon{\cat{W}}$ and $\Mon{\cat{W'}}$ are the
same as the objects of $\Mon{\cat{V}}$; that is, they are monads in
$\cat{B}$. But by using $\cat{W}$ or $\cat{W'}$ we obtain a more flexible
notion of a `map of monads' than we did in Example~(\ref{eg:strict-monads}):
a map in $\Mon{\cat{W}}$ is what Street called a \emph{monad functor} in
\cite{St}, and a map in $\Mon{\cat{W'}}$ is a \emph{monad opfunctor}.

\end{enumerate}

\section{Bimodules}	\label{sec:bim}

Bimodules have traditionally been discussed in the context of
bicategories. Thus given a bicategory $\cat{B}$, one constructs a new
bicategory $\Bim{\cat{B}}$ whose 1-cells are bimodules in $\cat{B}$ (see
e.g.\ \cite{CKW}). The drawback is that this is only possible when $\cat{B}$
has certain properties concerning the existence and behaviour of local
reflexive coequalizers.

Here we extend the $\fcat{Bim}$ construction from bicategories to
$\fc$-multicategories, which allows us to drop the technical assumptions. In
other words, we will construct an honest functor
\[
\fcat{Bim}: \fc\hyph\Multicat \go \fc\hyph\Multicat. 
\]
($\fc\hyph\Multicat$ is the category of (small) $\fc$-multicategories, with
maps defined in the obvious way.)

I would like to be able to, but at present cannot, place the $\fcat{Bim}$
construction in a more abstract setting: as it stands it is somewhat \emph{ad
hoc}. In particular, the definition does not appear to generalize to
$T$-multicategories for arbitrary $T$.

The theories of bimodules and enrichment interact in the following way:
given an $\fc$-multicategory $\cat{V}$, there is a canonically-defined functor
\[
\cat{V}\hyph\Cat \go \Bim{\cat{V}}\hyph\Cat.
\]
This is discussed at the end of the section, and provides lots of new
examples of enriched categories.

We first have to define $\fcat{Bim}$. Let $\cat{V}$ be an
$\fc$-multicategory: then the $\fc$-multicategory $\Bim{\cat{V}}$ is defined
as follows.
\begin{description}
\item[0-cells] A 0-cell of $\Bim{\cat{V}}$ is an $\fc$-multicategory map
$1\go \cat{V}$. That is, it is a 0-cell $x$ of $\cat{V}$ together with a
horizontal 1-cell $x\goby{t}x$ and 2-cells
\[
\begin{diagram}[height=2em]
x	&\rTo^{t}	&x			&\rTo^{t}	&x	\\
\dTo<{1}&		&\Downarrow\,\scriptstyle{\mu}&		&\dTo>{1}\\
x	&		&\rTo_{t}		&		&x	\\
\end{diagram}
\diagspace
\begin{diagram}[height=2em]
x	&\rEquals		&x	\\
\dTo<{1}&\Downarrow\,\scriptstyle{\eta}&\dTo>{1}\\
x	&\rTo_{t}		&x	\\
\end{diagram}
\]
satisfying the usual axioms for a monad, $\mu\of(\mu,1_t) =
\mu\of(1_t,\mu)$ and $\mu\of(\eta,1_t) = 1_t = \mu\of(1_t,\eta)$.

\item[Horizontal 1-cells]
A horizontal 1-cell $(x,t,\eta,\mu) \rTo (x',t',\eta',\mu')$ consists of a
horizontal 1-cell $x\goby{m}x'$ in $\cat{V}$ together with 2-cells
\[
\begin{diagram}[height=2em]
x	&\rTo^{t}	&x			&\rTo^{m}	&x'	\\
\dTo<{1}&		&\Downarrow\,\scriptstyle{\theta}&	&\dTo>{1}\\
x	&		&\rTo_{m}		&		&x'	\\
\end{diagram}
\diagspace
\begin{diagram}[height=2em]
x	&\rTo^{m}	&x'			&\rTo^{t'}	&x'	\\
\dTo<{1}&		&\Downarrow\,\scriptstyle{\theta'}&	&\dTo>{1}\\
x	&		&\rTo_{m}		&		&x'	\\
\end{diagram}
\]
satisfying the usual module axioms $\theta\of(\eta,1_m)=1_m$,
$\theta\of(\mu,1_m) = \theta\of(1_t,\theta)$, and dually for
$\theta'$, and the `commuting actions' axiom
$\theta'\of(\theta,1_{t'}) = \theta\of(1_t,\theta')$.

\item[Vertical 1-cells] A vertical 1-cell
$\vslob{(x,t,\eta,\mu)}{}{(\hat{x},\hat{t},\hat{\eta},\hat{\mu})} in
\Bim{\cat{V}}$ is a vertical 1-cell $\vslob{x}{f}{\hat{x}}$ in $\cat{V}$
together with a 2-cell
\begin{diagram}[height=2em]
x		&\rTo^{t}		&x		\\
\dTo<{f}	&\Downarrow\,\scriptstyle{\omega}&\dTo>{f}	\\
\hat{x}		&\rTo_{\hat{t}}		&\hat{x}	\\
\end{diagram}
satisfying $\omega\of\mu = \hat{\mu}\of(\omega,\omega)$ and a similar
equation for units. 

\item[2-cells]
A 2-cell 
\begin{diagram}[height=2em]
t_0	&\rTo^{m_1}	&t_1	&\rTo^{m_2}	&\ 	&\cdots	
&\ 	&\rTo^{m_{n}}	&t_n	\\
\dTo<{f}&		&	&		&\Downarrow	&	
&	&		&\dTo>{f'}\\
t	&		&	&		&\rTo_{m}	&	
&	&		&t'	\\
\end{diagram}
in $\Bim{\cat{V}}$, where $t$ stands for $(x,t,\eta,\mu)$, $m$ for
$(m,\theta,\theta')$, $f$ for $(f,\omega)$, and so on, consists of a
2-cell 
\begin{diagram}[height=2em]
x_0	&\rTo^{m_1}	&x_1	&\rTo^{m_2}	&\ 	&\cdots	
&\ 	&\rTo^{m_{n}}	&x_n	\\
\dTo<{f}&		&	&
&\Downarrow\,\scriptstyle{\alpha}& 	
&	&		&\dTo>{f'}\\
x	&		&	&		&\rTo_{m}	&	
&	&		&x'	\\
\end{diagram}
in $\cat{V}$, satisfying the `external equivariance' axioms
\begin{eqnarray*}
\alpha\of(\theta_1,1_{m_2},\ldots,1_{m_n}) 		&=&
\theta\of(\omega,\alpha)				\\
\alpha\of(1_{m_1},\ldots,1_{m_{n-1}},\theta'_n)		&=&
\theta'\of(\alpha,\omega')
\end{eqnarray*}
and the `internal equivariance' axioms
\begin{eqnarray*}
\lefteqn{\alpha\of(1_{m_1},\ldots,1_{m_{i-2}}, \theta'_{i-1}, 1_{m_{i}},
1_{m_{i+1}},\ldots,1_{m_n})		=}				\\
&&\alpha\of(1_{m_1},\ldots,1_{m_{i-2}}, 1_{m_{i-1}}, \theta_i,
1_{m_{i+1}},\ldots,1_{m_n}) 
\end{eqnarray*}
for $2\leq i\leq n$.

\item[Composition and identities]
For both 2-cells and vertical 1-cells in $\Bim{\cat{V}}$, composition is
defined directly from the composition in $\cat{V}$, and similarly identities.

\end{description}

Incidentally, the category formed by the objects and vertical 1-cells of
$\Bim{\cat{V}}$ is $\Mon{\cat{V}}$, the category of monads in $\cat{V}$
defined earlier. 

We have now defined an $\fc$-multicategory $\Bim{\cat{V}}$ for each
$\fc$-multicategory $\cat{V}$, and it is clear how to do the same thing for
maps of $\fc$-multicategories, so that we have a functor
\[
\fcat{Bim}: \fc\hyph\Multicat \go \fc\hyph\Multicat.
\]
Again, we have been rather eccentric in naming the `bimodules construction'
after what it does to the horizontal 1-cells rather than the objects: perhaps
we should call it the `monads construction'. We are, however, following the
traditional terminology.

\subsection*{Examples}

\begin{enumerate}

\item
Let $\cat{B}$ be a bicategory satisfying the conditions on local reflexive
coequalizers mentioned in the first paragraph of this section, so that it is
possible to construct a bicategory $\Bim{\cat{B}}$ in the traditional way. Let
$\cat{V}$ be the $\fc$-multicategory coming from $\cat{B}$. Then a 0-cell of
$\Bim{\cat{V}}$ is a monad in $\cat{B}$, a horizontal 1-cell $t\go t'$ is a
$(t',t)$-bimodule, and a 2-cell of the form
\begin{diagram}[height=2em]
t_0	&\rTo^{m_1}	&t_1	&\rTo^{m_2}	&\ 	&\cdots	
&\ 	&\rTo^{m_{n}}	&t_n	\\
\dTo<{1}&		&	&		&\Downarrow	&	
&	&		&\dTo>{1}\\
t_0	&		&	&		&\rTo_{m}	&	
&	&		&t_n	\\
\end{diagram}
is a map 
\[
m_n \otimes_{t_{n-1}}\cdots\otimes_{t_1} m_1 \go m
\] 
of $(t_n,t_0)$-bimodules, i.e.\ a 2-cell in $\Bim{\cat{B}}$. So if we discard
the non-identity 1-cells of $\Bim{\cat{V}}$, then the resulting
$\fc$-multicategory is precisely the $\fc$-multicategory associated with the
bicategory $\Bim{\cat{B}}$.

\item 
$\Bim{\Span} = \Prof$, where $\Span$ is the $\fc$-multicategory of sets,
functions, spans, etc, and $\Prof$ is the $\fc$-multicategory of categories,
functors, profunctors, etc (Examples~\ref{sec:fcms}(\ref{eg:Span})
and~(\ref{eg:Prof})). 

\item
$\Bim{\Sigma\Ab} = \Bimod$ (Example~\ref{sec:fcms}(\ref{eg:Bimod})). Here
$\Ab$ is regarded as a monoidal category under the usual tensor and
$\Sigma\Ab$ is as in Example~\ref{sec:fcms}(\ref{eg:mon-cat}); or
equivalently $\Ab$ is regarded as a multicategory with the usual multilinear
maps, and $\Sigma\Ab$ is as in Example~\ref{sec:fcms}(\ref{eg:multicat}). 

\item 
Similarly, $\Bim{\Sigma\Set} = \Action$
(Example~\ref{sec:fcms}(\ref{eg:Action})), with cartesian product giving the
monoidal category (or multicategory) structure on $\Set$.

\item 
It is possible to define an $\fc$-multicategory $\Span(\cat{E},T)$, for any
appropriate monad $T$ on a category $\cat{E}$, and then
$\Bim{\Span(\cat{E},T)}$ is the $\fc$-multicategory of $T$-multicategories
and maps, profunctors, etc, between them. See~\cite{GECM} or~\cite{OHDCT} for
details.

\end{enumerate}

We now show how the bimodules construction produces new enriched categories
from old. 

\begin{prop}
For any $\fc$-multicategory $\cat{V}$, there is a natural functor
\[
\tilde{(\:\:)}: \cat{V}\hyph\Cat \go \Bim{\cat{V}}\hyph\Cat, 
\]
preserving object-sets.
\end{prop}

\paragraph*{Proof}
Take a $\cat{V}$-enriched category $C$. We must define a
$\Bim{\cat{V}}$-enriched category $\tilde{C}$ with object-set $C_0$, and so,
for instance, we must define an object $\eend{\tilde{C}}{a}$ of
$\Bim{\cat{V}}$ for each $a\in C_0$. To do this we observe that $\eend{C}{a}$
has a natural monad structure on it: that is, we put
\[
\eend{\tilde{C}}{a} = (\eend{C}{a}, \ehom{C}{a}{a}, \ids_a, \comp_{a,a,a}).
\]
The rest of the construction is along similar lines; there is only one
sensible way to proceed, and it is left to the reader. (An abstract account
is in~\cite{OHDCT}).  
\hfill\ensuremath{\Box}

\subsection*{Examples}

\begin{enumerate}

\item 	\label{eg:Ab-to-Ring} 
Let $C$ be a category enriched (in the usual sense) in the monoidal category
$\Ab$ of abelian groups. Then the resulting $\Bimod$-enriched category
$\tilde{C}$ is as follows:
\begin{itemize}
\item $\tilde{C}_0$ is the set of objects of $C$
\item $\eend{\tilde{C}}{a}$ is the ring $\ehom{C}{a}{a}$, in which
multiplication is given by composition in $C$
\item $\ehom{\tilde{C}}{a}{b}$ is the abelian group $\ehom{C}{a}{b}$ acted on
by $\eend{\tilde{C}}{a} = \ehom{C}{a}{a}$ on the right and by
$\eend{\tilde{C}}{b} = \ehom{C}{b}{b}$ on the left, both actions being by
composition in $C$
\item composition and identities are as in $C$.
\end{itemize}

To illustrate the functoriality in the Proposition, take an $\Ab$-enriched
functor $F: C \go D$. This induces a $\Bimod$-enriched functor $\tilde{F}:
\tilde{C} \go \tilde{D}$ as follows:
\begin{itemize}
\item
$\tilde{F}: C_0 \go D_0$ is the object-function of $F$
\item
if $a\in C_0$ then $\tilde{F}_a$ is the ring homomorphism 
\[
\eend{\tilde{C}}{a} = \ehom{C}{a}{a} 
\go 
\ehom{D}{F(a)}{F(a)} = \eend{\tilde{D}}{\tilde{F}(a)}
\]
induced by $F$
\item
if $a,b\in C_0$ then
\[
\tilde{F}_{ab}: \ehom{\tilde{C}}{a}{b} = \ehom{C}{a}{b}
\go 
\ehom{D}{F(a)}{F(b)} = \ehom{\tilde{D}}{\tilde{F}(a)}{\tilde{F}(b)}
\]
is defined by the action of $F$ on morphisms $a\go b$.
\end{itemize}
Note that in general, the ring homomorphism $\tilde{F}_a$ is not the identity;
so the vertical 1-cells of $\Bimod$ get used in an essential way. This is
the reason why the Proposition does not hold if we work throughout with
bicategories rather than $\fc$-multicategories: $\tilde{(\:\:)}$ is defined on
objects of $\cat{V}\hyph\Cat$, but cannot sensibly be defined on morphisms. 

\item
The non-additive version of~(\ref{eg:Ab-to-Ring}) is that there is a
canonical functor
\[
\begin{array}{ccc}
\Cat	&\go		&\Action\hyph\Cat	\\
C	&\goesto	&\tilde{C}		
\end{array}
\]
which exists because, for instance, the set of endomorphisms on an object of
a category is naturally a monoid. 

\item
In the previous example, part of the construction was to take
$\eend{\tilde{C}}{a}$ to be the monoid of all endomorphisms of $a$ in
$C$. However, we could just as well take only the automorphisms of $a$, and
this would yield a different functor from $\Cat$ to $\Action\hyph\Cat$.

\item 
Applying the Proposition to $\cat{V}=\Span$ and recalling
Example~\ref{sec:enr}(\ref{eg:comma}), we obtain a functor
\[
(\ob\downarrow\Set) \go \Prof\hyph\Cat.
\]
What this does on objects is as follows. Take a category $D$, a set $I$, and
a function $\ob(D) \go I$. Then in the resulting $\Prof$-enriched category $E$
we have $E_0=I$; $\eend{E}{i}$ is the full subcategory of $D$ whose objects are
those lying over $i\in I$; and $\ehom{E}{i}{j}$ is the profunctor
\[
\begin{array}{ccc}
\eend{E}{i}^{\op} \times \eend{E}{j}	&\go		&\Set	\\
(d,d')					&\goesto	&D(d,d').
\end{array}
\]
Composition and identities are defined as in $D$.

\item
To get more examples of $\Prof$-enriched categories we can modify the
previous example, taking $\eend{E}{i}$ to be \emph{any} subcategory of $D$
whose objects are all in the fibre over $i$. Here are two specific
instances (each with a vague flavour of topological quantum field theory
about them). In the first, $E_0$ is the set $\mathbb{N}$ of natural
numbers, $\eend{E}{n}$ is the category of $n$-dimensional Hilbert spaces
($=$ complex inner product spaces) and isometries, and $\ehom{E}{m}{n}$
sends $(H,H')$ to the set of all linear maps $H \go H'$. In the second,
$E_0=\mathbb{N}$ again, and we replace Hilbert spaces by differentiable
manifolds, isometries by diffeomorphisms, and linear maps by differentiable
maps.

\end{enumerate}

\subsection*{Acknowledgements}

I am grateful to Martin Hyland for his many words of wisdom, and for
encouraging me to develop these ideas. I would also like to thank the
organizers of CT99 at Coimbra for the opportunity to present this work and to
discuss it with others. I thank the referee for a number of useful
suggestions. 

All but one of the diagrams in this document were generated using Paul
Taylor's commutative diagrams package.


\end{document}